\documentclass[12 pt]{article} 
\usepackage{amsmath,epsf}
\usepackage[latin1]{inputenc}
\usepackage{amsfonts}

\title{Monomial bases related to the $n!$ conjecture}
\author{Jean-Christophe Aval\\ Laboratoire A2X\\ Universit\'e Bordeaux 1\\ 351 cours de la Lib\'eration\\ 33405 Talence cedex\\ e-mail : {\tt aval@math.u-bordeaux.fr}}
\date{}

\font \petit=cmr9

\def\N{{\mathbb N}}

\def\QED{{\flushright $\Box$\\}}
\def\proof{{\it Proof. }}

\numberwithin{equation}{section}

\begin{document}

\maketitle


{\bf Abstract:} {\petit The purpose of this paper is to find a new way to prove the $n!$ conjecture for particular partitions. The idea is to construct a monomial and explicit basis for the space $M_{\mu}$. We succeed completely for hook-shaped partitions, i.e., $\mu=(K+1,1^L)$. We are able to exhibit a basis and to verify that its cardinality is indeed $n!$, that it is linearly independent and that it spans $M_{\mu}$. We derive from this study an explicit and simple basis for $I_{\mu}$, the annihilator ideal of $\Delta_{\mu}$. This method is also successful for giving directly a basis for the homogeneous subspace of $M_{\mu}$ consisting of elements of $0$ $x$-degree.} 

\section{Introduction}

Let $\mu=(\mu_1\ge\mu_2\ge\dots\ge\mu_k>0)$ be a partition of $n$. We shall identify $\mu$ with its Ferrers diagram (using the French notation). To each cell $s$ of the Ferrers diagram, we associate its coordinates $(i,j)$, where $i$ is the height of $s$ and $j$ the position of $s$ in its row. The pairs $(i-1,j-1)$ occurring while $s$ describes $\mu$ will be briefly referred to as the set of the biexponents of $\mu$. Now let $(p_1,q_1),\ldots,(p_n,q_n)$ denote the set of biexponents arranged in lexicographic order and set 
$$\Delta_{\mu}(x,y)=\Delta_{\mu}(x_1,\ldots,x_n;y_1,\ldots,y_n)=\det(x_i^{p_j}y_i^{q_j})_{i,j=1\ldots n}.$$
Let $M_{\mu}$ be the collection of polynomials in the variables $x_1,\ldots,$ $x_n;$ $y_1,\ldots,y_n$ obtained by taking the linear span of all the partial derivatives of $\Delta_{\mu}$. Formally we may write 
$$M_{\mu}={\cal L}\{\partial_x^a\partial_y^b\Delta_{\mu}(x,y);\ a,b\in \N^n\}$$
where $\partial_x^a=\partial_{x_1}^{a_1}\ldots\partial_{x_n}^{a_n}$ and $\partial_y^b=\partial_{y_1}^{b_1}\ldots\partial_{y_n}^{b_n}$.
Then the $n!$ conjecture can be stated as follows.

\noindent
{\bf Conjecture 1 ($n!$ conjecture):} {\it Let $\mu$ be a partition of $n$, then  $\dim M_{\mu}=n!$}.

This conjecture, stated by A. Garsia and M. Haiman is central for their study of Macdonald polynomials (cf. [5], [6]). To be more precise, Macdonald introduced in [12] a new symmetric function basis and associated Macdonald-Kostka coefficients $K_{\lambda\mu}(q,t)$, which are a priori rational functions in $q,t$. Macdonald conjectured that:

\noindent
{\bf Conjecture 2 (MPK conjecture):} {\it The functions $K_{\lambda\mu}(q,t)$ are polynomials with non-negative integer coefficients.}

Looking for a representation theoretical setting for the Macdonald basis, A. Garsia and M. Haiman made the following conjecture:

\noindent
{\bf Conjecture 3 ($C=\tilde H$ conjecture):} {\it For the diagonal action of $S_n$, $M_{\mu}$ is a bigraded version of the left regular representation. Moreover, if $C_{\lambda\mu}(q,t)$ denotes the bigraded multiplicity of the character $\chi_{\lambda}$ in the bigraded character of the module $M_{\mu}$ then: $C_{\lambda\mu}(q,t)=K_{\lambda\mu}(q,1/t)t^{n(\mu)}$, where $n(\mu)=\sum_{i=1}^k(i-1)\mu_i$.}

Conjecture 3 clearly implies Conjecture 1 and 2. M. Haiman [9] using Hilbert schemes theory recently proved that the $n!$ conjecture actually implies the $C=\tilde H$ conjecture. A part of the MPK conjecture is that the $K_{\lambda\mu}(q,t)$ are polynomials, which is not obvious from their definition. This part has been recently proved in several independent papers (cf. [7], [8], [10], [11], [14]).

When $\mu=(1^n)$ or $\mu=(n)$, $\Delta_{\mu}$ reduces to the Vandermonde determinant in $x$ and $y$ respectively. In these cases, it is a classical result (see [2]) that $\dim M_{\mu}=n!$. But although this conjecture has been verified by computer for small partitions up to $n=8$ and proved for some special cases (cf. [1], [4], [6], [13]), it has not been established in full generality. Several methods have been developed to prove the $n!$ conjecture but none of them has been able to give a proof in more than some special cases. 

In this paper our goal is to propose a new method to prove the $n!$ conjecture for some particular partitions. We want to construct explicit bases for the space $M_{\mu}$. These bases are made of monomial derivatives of $\Delta_{\mu}$. We present here how we are able to do it for hook-shaped partitions, i.e., $\mu=(K+1,1^L)$ with $K+L+1=n$. In section 2 we describe the way to construct the basis and prove that its cardinality is $n!$. In the third section we show that our family spans $M_{\mu}$. Moreover, we derive from that proof an explicit and simple basis for $I_{\mu}$, the annihilator ideal of $\Delta_{\mu}$. In the fourth section we prove by a completely new method that the elements of our basis are linearly independent. In section 5 we explain how this method is also successful for the homogeneous subspace of $M_{\mu}$ consisting of elements of $0$ $x$-degree. We obtain in fact a direct way to construct a basis for this subspace. 

\section{Construction and enumeration}

Let $\mu$ be a partition of $n$ whose Ferrers diagram is a hook, i.e., $\mu=(K+1,1^L)$ with $K+L+1=n$.

\subsection{Construction}

Let us take an horizontal axis.  A ``shape'' associated to $\mu$ is constructed the following way: suppose the line has room for $K+L$ spaces. Choose $K$ of these spaces to be $y$-columns and $L$ to be $x$-columns. In the $y$-columns place stacks of boxes above the line of height $K, K-1,\ldots, 1$ arranged in decreasing order. In the $x$-columns place stacks of boxes of decreasing depth $L, L-1,\ldots , 1$ below the line. 

Here is an example of shape: 
\vskip 0.5 cm

\centerline{
\epsffile{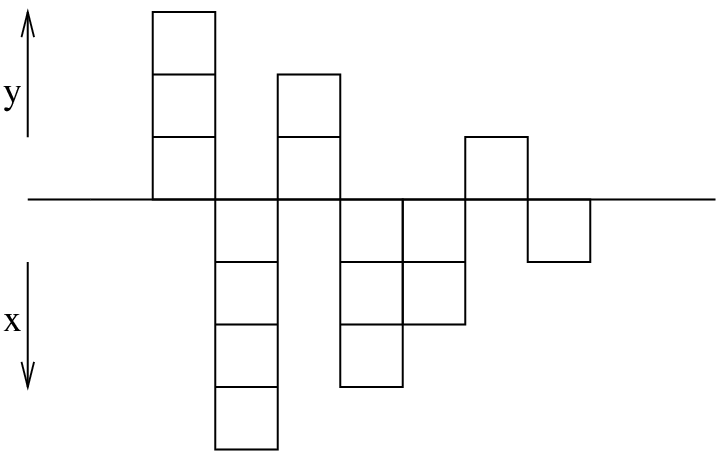}}

\vskip 0.3 cm

associated to the partition:

\centerline{
\epsffile{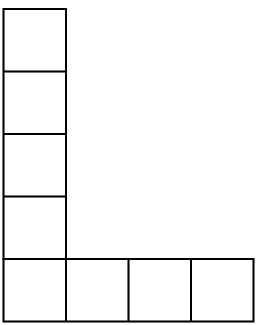}}

\vskip 0.3 cm

We shall now put crosses in the cells of the shape to obtain ``drawings''. As we shall not distinguish two drawings with the same number of crosses in each column, we put the crosses near the axis. The rules for putting crosses in a drawing are the following:
\begin{enumerate} 
\item the number of crosses in the $x$-columns is any number (not greater than the depth of the column);
\item the number of crosses in the $y$-columns depends on the $x$-crosses. For a column which has no $x$-column to its right, the number of crosses is not greater than the height of the column. In the other case, we look at the first ``plain'' $x$-column on the right; i.e., the first column which has only crosses (full $x$-column) or only white cells (empty $x$-column). There is always one, at least the $x$-column of depth one. Then:
\begin{itemize}
\item if it is all white, then we impose at least one cross in the $y$-column. 
\item if it is all crossed, then we impose at least one white cell in the $y$-column. 
\end{itemize}
\end{enumerate}

{\bf Remark 1:} The family of drawings that we defined is invariant under the operator that inverts the white cells and the crosses. We call this operator flip (it is different from the flip introduced by A. Garsia and M. Haiman in [6], that we denote from now on by Flip).

Here we give an example of drawing with crosses: 
\vskip 0.5 cm

\centerline{
\epsffile{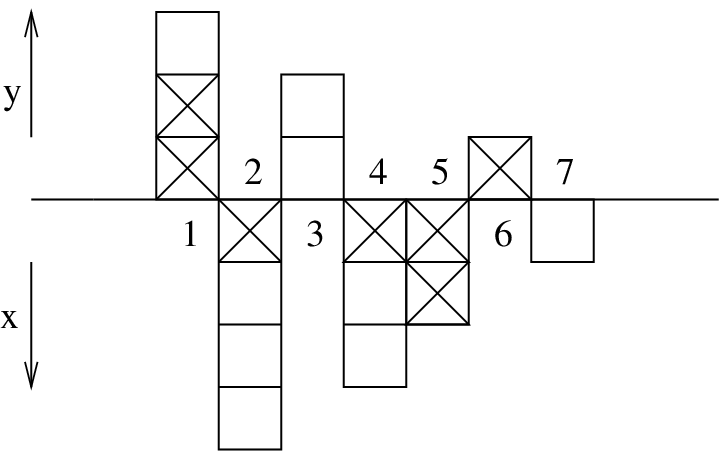}}

\vskip 0.3 cm

Once we have defined the drawings (with crosses), we define associated derivative operators. We give an index to the places of the drawing from left to right and from $1$ to $n-1$. Then to each $x$-cross in place $i$, we derive once with respect to $x_i$. We do the same thing for the $y$-crosses. For example, for the last drawing, the associated derivative operator is: $\partial_D= \partial y_1^2\partial x_2 \partial x_4\partial x_5^2\partial y_6$.

\subsection{Enumeration}

We shall denote by ${\cal D}$ the set of drawings that we defined in the previous subsection. We now verify that its cardinality is $n!$.

As the number of choices for the $y$-columns depends only on the shape of the drawing (and not on the $x$-crosses), we can write that the cardinality equals  the following expression, where $k_1$ denotes the number of $y$-columns on the right of the last $x$-column:
$$\sum_{k_1+k_2=K} 2\cdot3\cdots(k_1+1)\cdot(k_1+1)\cdots(k_1+k_2)\cdot(L+1)!\ {{k_2+L-1} \choose {k_2}}$$
$$=L(L+1)K! \sum_{k_2=0}^K \frac {(k_2+L-1)!} {k_2!} (K+1-k_2)$$
$$=(L+1)!K!\sum_{k_2=0}^K {{L-1+k_2}\choose{L-1}} {{K+1-k_2}\choose{1}}$$
$$=(L+1)!K!{{K+L+1}\choose{L+1}}=(K+L+1)!$$
by the Chu-Vandermonde formula ([3], p. 163).

\section{Proof that the family spans $M_{\mu}$}

We show here that $\{\partial_D \Delta_{\mu}\}_{D\in{\cal D}}$ spans $M_{\mu}$. We begin by studying $I_{\mu}$, the annihilator ideal of $\Delta_{\mu}$.

\subsection{Study of $I_{\mu}$}

For $P,\ Q$ two polynomials, we write $P\equiv Q$ if $P(\partial)\Delta_{\mu}=Q(\partial)\Delta_{\mu}$, i.e., $P-Q\in I_{\mu}$ ($P(\partial)$ corresponds to the substitution: $x_i\rightarrow \partial x_i,\ y_i\rightarrow \partial y_i$). We denote as usual by $h_k$ the $k$-th complete homogeneous symmetric function. Let also $X$ denote a subset of $(x_1,x_2,\ldots,x_n)$, $Y$ a subset of $(y_1,y_2,\ldots,y_n)$, $|X|$ and $|Y|$ their cardinality. We also set $\bar{X}=\prod_{x\in X}x$ and $\bar{Y}=\prod_{y\in Y}y$.

We first notice that:
\begin{enumerate}
\item for all $1\le i\le n$, $x_iy_i\equiv 0$;
\item $\bar{X}\equiv 0$ as soon as $|X|>L$;
\item $\bar{Y}\equiv 0$ as soon as $|Y|>K$;
\item for any symmetric homogeneous polynomial $P$ of positive degree, $P\equiv 0$.
\end{enumerate}
The fourth relation is well known (cf. [2]). The others are clear by observing the elements in the determinantal form of $\Delta_{\mu}$ when $\mu=(K+1,1^L)$.

\vskip 0.2 cm
{\bf Proposition 1:} 
{\it $$h_{k}(Y)\equiv 0$$
as soon as $k>0$ and $k+|Y|>n$.}

\proof
It is easily proved by an induction based on $h_k(y_1,\ldots,y_{n})\equiv 0$ for all $k>0$. We have indeed $h_{1}(Y_n)\equiv 0$, where $Y_n=(y_1,\ldots,y_n)$ and for any $y\not\in Y$:
$$h_k(Y,y)=h_k(Y)+yh_{k-1}(Y,y).$$
\QED

\vskip 0.2 cm
{\bf Proposition 2:} 
{\it $$\bar{Y}h_{k}(Y')\equiv 0$$
as soon as $k>0$, $k+|Y|>K$ and $Y\subset Y'$.}

\proof
Proposition 2 is proved by decreasing induction on $|Y'|$.

We observe that the result is true for $|Y'|$ equal to $K+1$ and $K$. Let $Y$ and $Y'$ satisfy the hypotheses and assume the result is true down to $|Y'|+1$. We write for all $y_i \not \in Y'$:
$$h_{k}(Y',y_i)\equiv h_{k}(Y')+y_{i}h_{k-1}(Y',y_{i}),$$
thus, by induction if $k>1$ we obtain the following relation; this relation is obvious if $k=1$ because this implies that $|Y|\ge K$:
$$\bar{Y}h_{k}(Y',y_{i})\equiv \bar{Y}h_{k}(Y').$$
Once we have this relation the conclusion easily follows by an increasing induction on $|Y'|$ (for example up to $n$).
\QED

\vskip 0.2 cm
{\bf Proposition 3:} 
{\it $$h_{k}(Y)h_{l}(X)\equiv 0$$
as soon as $k>0$, $l>0$, $k+l+|Y|+|X|\ge 2n$ and $X\subset Y$ or $Y\subset X$.}

\proof 
We only show the result when $k+|Y|=n$ and $l+|X|=n$ (the other cases are consequences of Proposition 1).

It is in fact proved as Proposition 1 by a simple induction based on:
$$h_1(x_1,\ldots,x_{n-1}) h_1(y_1,\ldots,y_{n-1})\equiv 0$$
which is a consequence of Proposition 1 and $x_ny_n\equiv 0$.
\QED

\vskip 0.2 cm
{\bf Proposition 4:} 
{\it $$h_{k}(Y)h_{l}(X)\equiv 0$$
as soon as $k>0$, $l>0$ and 
\begin{itemize}
\item either $Y\subset X$ and $k+l+|Y|>n$,
\item or  $X\subset Y$ and $k+l+|X|>n$.
\end{itemize}
}

\proof
This is proved by induction on $\alpha=2n-(k+|Y|+l+|X|)$.

The case $\alpha\le 0$ reduces to Proposition 3.

Suppose the result is true up to $\alpha-1$ and $2n-(k+|Y|+l+|X|)=\alpha>0$. By symmetry, we shall assume that $Y\subset X$ and $k+l+|Y|>n$. If $l>1$, then for any $x_i\not \in X$, we write: 

$$h_{k}(Y)h_{l}(X)\equiv h_{k}(Y)h_{l}(X,x_{i})-x_{i}h_{k}(Y)h_{l-1}(X,x_i)$$
$$\equiv h_{k}(Y)h_{l}(X,x_{i})-x_{i}h_{k}(Y,y_{i})h_{l-1}(X,x_{i})\equiv 0$$
by induction.

If $l=1$, then $|Y|+k\ge n$ and we write for any $x_i\not\in X$:
$$h_{k}(Y)h_{1}(X)\equiv h_{k}(Y,y_i)h_{1}(X)-y_{i}h_{k-1}(Y)h_{1}(X,x_i).$$
The first term is zero by Proposition 1. The second term is proved to be also zero by increasing induction on $|X|$ (up to $n$), since $n-k\le |Y|\le |X|\Rightarrow n-|X|\le k$.
\QED

\subsection{Application}

We shall show here that any monomial derivative of $\Delta_{\mu}$ is a linear combination of the derivatives: $\{\partial_D\Delta_{\mu}\}_{D \in {\cal D}}$ (derivatives corresponding to drawings, i.e., the family defined in section 2).

\vskip 0.3 cm

{\bf Theorem 1: }{\it $\{\partial_D\Delta_{\mu}\}_{D \in {\cal D}}$ spans $M_{\mu}$.}a

\vskip 0.3 cm

\proof
It is clear that any monomial can be associated to a diagram of crosses (by the same process as in paragraph 2.1), and let $D$ be such a diagram which is not a drawing.

We look at the rightmost ``anomaly'', that is the rightmost place where the diagram $D$ associated to the monomial can not be a drawing (we call this place ``guilty'').

\begin{description}
\item{-} Case 1: the diagram $D$ could not be put in a set of ordered columns (i.e., in the shape of a drawing). This case gives four subcases. Assume the guilty column is a $y$-column. We can not put another $y$-column on the right. Either because each $y$-column on the left has a cross (case 1a), or because there is no first plain and white $x$-column on the right (case 1b). If the guilty column is an $x$-column, we are led to cases 1c (each $x$-column on the left has a cross) and 1d (there is no first plain and white $y$-column on the right). Since the rules are not involved here the problems are symmetric for $x$ or $y$.
\item{-} Case 2: the diagram $D$ could be put in a set of ordered columns but the rules are broken. Either for the white cells (case 2a), or for the crosses (case 2b).
\end{description}

We shall prove, using the propositions of the last section, that the monomial associated to the diagram $D$ can be written modulo $I_{\mu}$ as a linear combination of monomials strictly smaller with respect to the lexicographic order ($x_1<x_2<\dots<x_n<y_1<\dots<y_n$). We look at each case that we have mentioned above. 

\begin{itemize}
\item Case 1b with no $x$-column on the right is solved by Proposition 1, as well as case 1d with no $y$-column on the right.

\item Cases 1a and 1c are symmetric and treated by Proposition 2: we note that the height of the $h$-th $y$-column is $K-h+1$. If it has $k+1$ crosses, there is a problem if $k+h>K$. It then can be treated by Proposition 2: we take $Y'=Y=\{i_1<\ldots<i_h\}$, to be the places of the first $h$ $y$-columns, each of which has at least one cross. The monomial is a multiple of 
$$\bar{Y}y_{i_h}^k\equiv \bar{Y}(y_{i_h}^k-h_k(Y))$$
and all monomials in the expansion of the right side are lexicographically smaller than the monomial on the left side. 

\item Case 2a is immediately settled by inverting the involved columns. 
\end{itemize}

Therefore the only remaining cases are case 1b (resp. 1d) with a first full $x$- (resp. $y$-) column on the right and case 2b.

\begin{itemize}
\item Let us first study the case 2b.

\vskip 0.5 cm

\centerline{
\epsffile{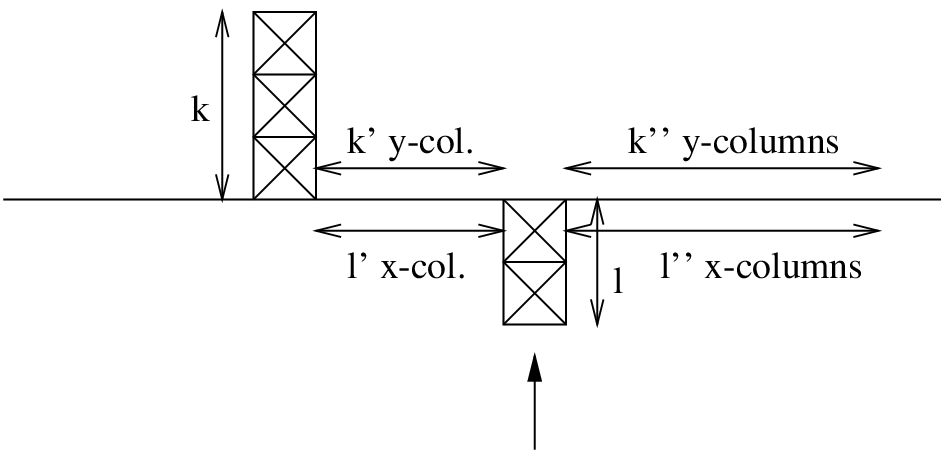}}

\vskip 0.3 cm

We observe that there is a problem if one has simultaneously:
\begin{itemize}
\item $k=k'+k''+1$,
\item $l=l''+1$,
\item there is a cross in each of the $l'$ $x$-columns between the two columns appearing on the figure.
\end{itemize}

Let:
\begin{itemize}
\item $Y$ denote the places on the left of the $y$-column plus the place of the $y$-column plus the $l'$ places of the $x$-columns between the $y$- and the $x$-column on $D$ plus the place of the $x$-column,
\item $X$ denote the places on the left of the $x$-column plus the place of the $x$-column itself,
\item $X'$ denote the places of the $l'$ $x$-columns between the $y$- and the $x$-column of $D$.
\end{itemize}
We shall be able to express the monomial corresponding to this $D$ as a linear combination of monomials stricly smaller with respect to the lexicographic order if we establish that
$$h_{k}(Y)h_{l}(X)\equiv 0.$$
Indeed the leading monomial of $\bar X'h_{k}(Y)h_{l}(X)$ (for the lexicographic order), in which we delete the multiples of $x_iy_i$ for any $i$, is a divisor of the monomial associated to $D$. 

We want to apply Proposition 4 with $|Y|=n-(k'+k''+l''+1)$ and $|X|=n-(k''+l''+1)$. We have $Y\subset X$ and we calculate:
$$k+l+|Y|-n=1>0.$$
Hence we are done in this case.

\item Let us now consider the case 1d with a first full $y$-column.

\vskip 0.5 cm

\centerline{
\epsffile{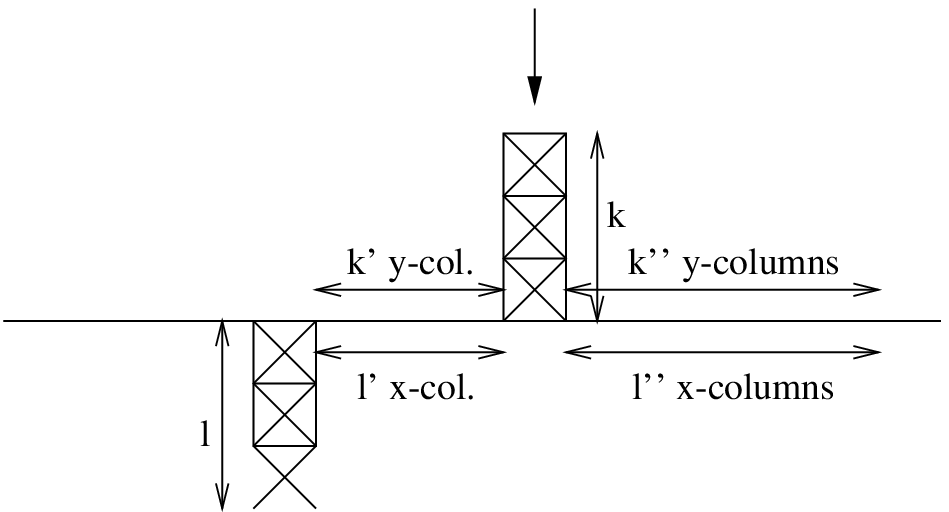}}

\vskip 0.3 cm

Here a problem occurs if:
\begin{itemize}
\item $k=k''+1$,
\item $l\ge l'+l''+2$,
\item there is a cross in each of the $k'$ $y$-columns between the two columns appearing on the figure.
\end{itemize}

We proceed as in the previous case. We want to use Proposition 4 to show that
$$h_{k}(Y)h_{l}(X)\equiv 0$$
with $Y$ corresponding to all the places strictly left of the $y$-column on the diagram $D$ and $X$ corresponding to all the places up to the $x$-column, plus the places of the $k'$ $y$-columns between the $x$-column and the $y$-column. 

We want to apply Proposition 4 with $|X|=n-(l'+k''+l''+2)$ and $|Y|=n-(k''+l''+1)$. We have $X\subset Y$ and we compute:
$$k+l+|X|-n\ge 1.$$
Thus this case is also settled.

\item It remains to observe that the case 1b with a first full $x$-column is treated by case 2b.
\end{itemize}

The proof of Theorem 1 is now complete.
\QED

\subsection{Conclusion}

We can deduce from what precedes a basis for the ideal $I_{\mu}$ when $\mu$ is a hook, since the first relations exposed at the beginning of the study of $I_{\mu}$ were sufficient to prove that our family is a basis of $M_{\mu}$.

{\bf Theorem 2:} {\it If we denote by $\langle G\rangle$ the ideal generated by a set $G$, then for $\mu$ a hook partition of $n$, we have:
$$I_{\mu}=\langle h_i(X_n),\ 1\le i\le n;\ h_i(Y_n),\ 1\le i\le n;$$ 
$$x_iy_i,\ 1\le i\le n;\ \bar{X},\ |X|=L+1;\ \bar{Y},\ |Y|=K+1\rangle.$$}

\proof
To prove this we assume that the previous ideal (we denote it by $I$) is not equal to $I_{\mu}$, so that there is a polynomial $P$ in $I_{\mu}\backslash I$. According to the proof of Theorem 1, we can decompose it as $P=A+Q$, where $A$ is a linear combination of monomials of our family and $Q$ is an element of $I$. Taking the derivatives and applying it to $\Delta_{\mu}$, we obtain $A(\partial)\Delta_{\mu}=0$. As we shall see in section 4, this implies $A=0$, and $P=Q\in I$.
\QED

\section{Proof of the independence}

\subsection{Exposition and reduction of the problem}

We shall now prove that our family is an independent set.

Since the derivative operator associated to a drawing $D$ depends only on the crosses and not on the shape of the drawing, we define $S$ as the diagram consisting only of the crosses of $D$. We also define $T$ as the diagram consisting of the white cells (a ``complement'' of $S$).

Let ${\cal S}$ denote the set of $S$ when $D$ varies in ${\cal D}$, the set of drawings defined in section 2.

For example, for the drawing in section 2, we have:
\vskip 0.5 cm

\centerline{
\epsffile{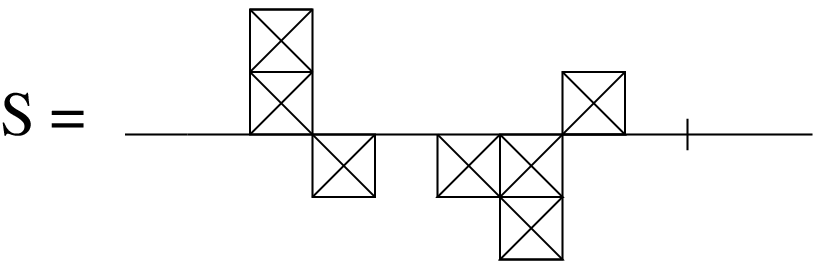}}

\vskip 0.3 cm
\vskip 0.5 cm

\centerline{
\epsffile{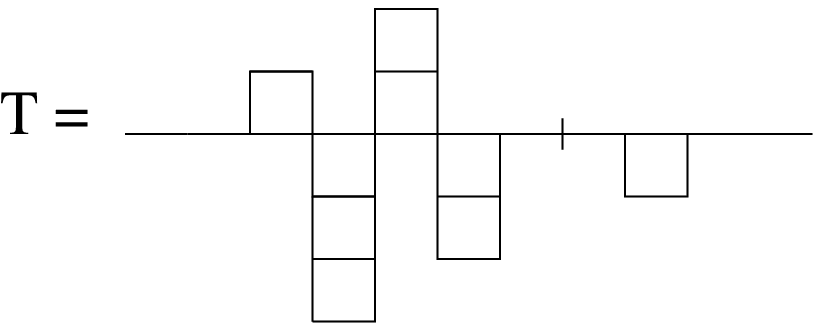}}

\vskip 0.3 cm

Let now $\partial_S$ and $\partial_T$ denote the derivative operators associated to $S$ and $T$ (after putting crosses in all the cells of $T$). 

\vskip 0.3 cm
{\bf Theorem 3:} 

{\it The family $\{\partial_S.\Delta_{\mu}\}_{S \in {\cal S}}$ is an independent set.}

\vskip 0.5 cm
{\bf Lemma 1:}

{\it
$S$ or $T$ determines the drawing from which it comes.}

\proof
Indeed, we can recontruct the shape of the drawing from $S$ by proceeding from left to right. The method is the following: if there are crosses at the place we are looking at, we complete the column with respect to the size of the successive columns. If there is no cross, we look at the $x$-crosses on the right: if they can fit in with one $x$-column missing, then we put an $x$-column at the empty place, else we put a $y$-column.

The method is the same for $T$ since the family is invariant under flip.
\QED

\vskip 0.3 cm
Let us now show that the family is linearly independent. Let us begin with some definitions. Let $D=(S,T)$ and $D_1=(S_1,T_1)$ be two different drawings; we shall say that $D_1$ is a son of $D$ if $\partial_T\circ\partial_{S_1}.\Delta_{\mu}\in {\mathbb Z}\backslash\{0\}$. We shall denote by $T+S_1$ the figure corresponding to the superposition (place by place) of the cells of $T$ and $S_1$ (all these cells being crossed). If we repeat this process, we obtain the notion of descendant.

\vskip 0.3 cm
{\bf Lemma 2:}

{\it To show the independence, it is sufficient to prove that a drawing can not be its own descendant (i.e., there is no ``loop'').}

\proof
We assume we have a relation of dependence: $\sum_{S} c_S\partial_S.\Delta_{\mu}=0$, that the coefficients are not all zero, and that there is no loop. Then we take a $S_0$ for which $c_{S_0}\neq 0$. If $S_0$ has no son or if they have all $c_S$ equal to zero, we obtain a contradiction by applying $\partial_{T_0}$ to the relation and by looking at the constant term of the result. If $S_0$ has a son $S_1$ for which $c_{S_1}\neq 0$, we repeat with $S_1$. As the set is finite and there is no loop, we certainly obtain a $S'$ which gives a contradiction.
\QED

\vskip 0.2 cm

So we have to prove that there is no loop. It is sufficient to show that a drawing $D=(S,T)$ is different from all its descendants that have the same shape (i.e., the $x$-columns at the same places). Let $D'=(S',T')$ be a descendant of $D$ that has the same shape. We want to show that $D\neq D'$.

\subsection{Definition of completeness}

To explain this notion, let $D_1$ denotes a drawing and $D_2$ one of its sons. We define on the places of $D_2$ a notion of ``completeness'' (relative to  $D_1$ too) as follows:
We say that the first $k$ places of $D_2$ are complete if the heights of the $y$-columns of $T_1+S_2$ in these $k$ places and read from left to right are $K,K-1,K-2,\ldots$ and if we have the same for $x$-columns.

We want now to obtain a (more quantitative) characterization of the completeness. To do this we need to  introduce some more definitions.

We look at the left parts (made of the first $k-1$ places) of $D_1$ and $D_2$. We define $d$ as the difference between the number of times where a $y$-column of $D_1$ has been replaced in $D_2$ by a white $x$-column and the number of times where an $x$-column of $D_1$ has been replaced in $D_2$ by a white $y$-column. We also define $d'$ as the difference between the number of times where a crossed $y$-column of $D_1$ has been replaced in $D_2$ by an $x$-column and the number of times where a crossed $x$-column of $D_1$ has been replaced in $D_2$ by a $y$-column. We should note that $d$ and $d'$ are relative to $k-1$.

Since the problem is symmetric with respect to $x$ and $y$ (as long as we do not use the rules of construction), we shall only examine the case where we derive with respect to $y_k$, i.e., where there is a $y$-column at the $k$-th place of $T_1+S_2$. The symmetric case has a similar characterization (with opposite signs for $d$ and $d'$). We now introduce the following notations: $b_1$ (resp. $b_2$) denotes the number of white cells at place $k$ in $D_1$ (resp. $D_2$) and $c_1$ (resp. $c_2$) the number of crosses. The characterization can now be stated as follows:

{\bf Characterization :}

{\it
If the first $k-1$ places are complete, the $k$-th is complete if one of the following conditions is verified:
\begin{enumerate}
\item at place $k$ in $D_1$ and $D_2$ there is a $y$-column and $b_2=b_1+d$ and $c_2=c_1+d'$ (each of these equalities easily implies the other);
\item at place $k$, there is a crossed $x$-column in $D_1$ (i.e., $b_1=0$) and a $y$-column in $D_2$, and $b_2=d$;
\item at place $k$, there is a $y$-column in $D_1$ and a white $x$-column in $D_2$ ($c_2=0$), and $c_1=-d'$.
\end{enumerate}}

\proof
To prove this result, we begin by observing that we can not have $x$- and $y$- cells at the same place in $T_1+S_2$: when $\mu$ is a hook, we have $\partial x_i\partial y_i\Delta_{\mu}=0$. There are in fact three possibilities for the columns at place $k$:
\begin{enumerate}
\item $D_1$ and $D_2$ have a $y$-column;
\item $D_1$ has a crossed $x$-column and $D_2$ a $y$-column;
\item $D_1$ has a $y$-column and $D_2$ a white $x$-column.
\end{enumerate}
We deal with these three cases.
\begin{enumerate}
\item Case 1: if in $T_1+S_2$ the heights of the $y$-columns in the first $k-1$ places are $K$, $K-1$, $\ldots$, $K-l+1$ and if our $y$-column is the $h$-th of $D_2$, we observe that $l=h-1+d$. The height of the $y$-column of $T_1+S_2$ at place $k$ is at most $K-l$. But if we observe that the height of the $h$-th $y$-column of $D_2$ is $K-h+1$, we obtain:
$$b_1+c_2\le K-l=K-h+1-d=b_2+c_2-d.$$
Hence $b_2 \ge b_1+d$ and equality holds when it is complete. As $b_2+c_2=b_1+c_1+d+d'$, the equality $c_2=c_1+d'$ holds too.
\item Case 2: this case is treated like Case 1.
\item Case 3: the reasoning is similar to Case 1. If our $y$-column is the $h$-th of $D_1$ and if in the first $k-1$ places of $T_1+S_2$ the successive $y$-columns have height $K$, $K-1$, $\ldots$,  $K-l+1$, then $l=h-1-d'$. As the height of the $y$-column at place $k$ of $T_1+S_2$ is at most $K-l$, we deduce that $c_1\ge -d'$, with equality corresponding to completion.\QED
\end{enumerate}

\vskip 0.2 cm
{\bf Remark 2:} If the first $k-1$ places are complete but not the $k$-th, we observe easily that it corresponds to an increasing of the number of white cells in $D_2$. We have indeed seen in the proof that $b_2\ge b_1+d$ and $c_2\le c_1+d'$ in Case 1 and similar inequalities in Cases 2 and 3. 

\vskip 0.2 cm
{\bf Remark 3:} We observe that the Cases 2 and 3 can not happen simultaneously since we can not have at the same place a crossed column in $D_1$ and a white column in $D_2$ (there is at least one cell at each place).

Once we have obtained this characterization of completeness, we shall use it to progress in the proof of Theorem 3.

\subsection{Application}

\vskip 0.3 cm
{\bf Lemma 3:}

{\it If we have completeness on the first $k$ places along the chain between two drawings $D$ and $D'$ with the same shape, then the sum of the $d$ along the chain is equal to zero, as well as the sum of the $d'$ ($d$ and $d'$ relative to the first $k$ places).}

We will first apply this result in the following lemma and prove it after Lemma 4.

\vskip 0.3 cm
{\bf Lemma 4:}

{\it If we have completeness on the first $k$ places between $D$ and $D'$, then these two drawings are identical on the first $k$ places.}

\proof
To prove this result we shall use Lemma 3. 
Indeed we notice that if we keep either an $x$- or a $y$-column at place $k$ along the chain between $D$ and $D'$, the result is obvious since (by Lemma 3) the sum of the $d$ is equal to zero. With natural notations, we have:  $b'=b+\sum d=b$. Now, if the ``shape'' of the column at place $k$ changes, let us observe the two following cases (by symmetry we look at the changes for a $y$-column):
 
\centerline{
\epsffile{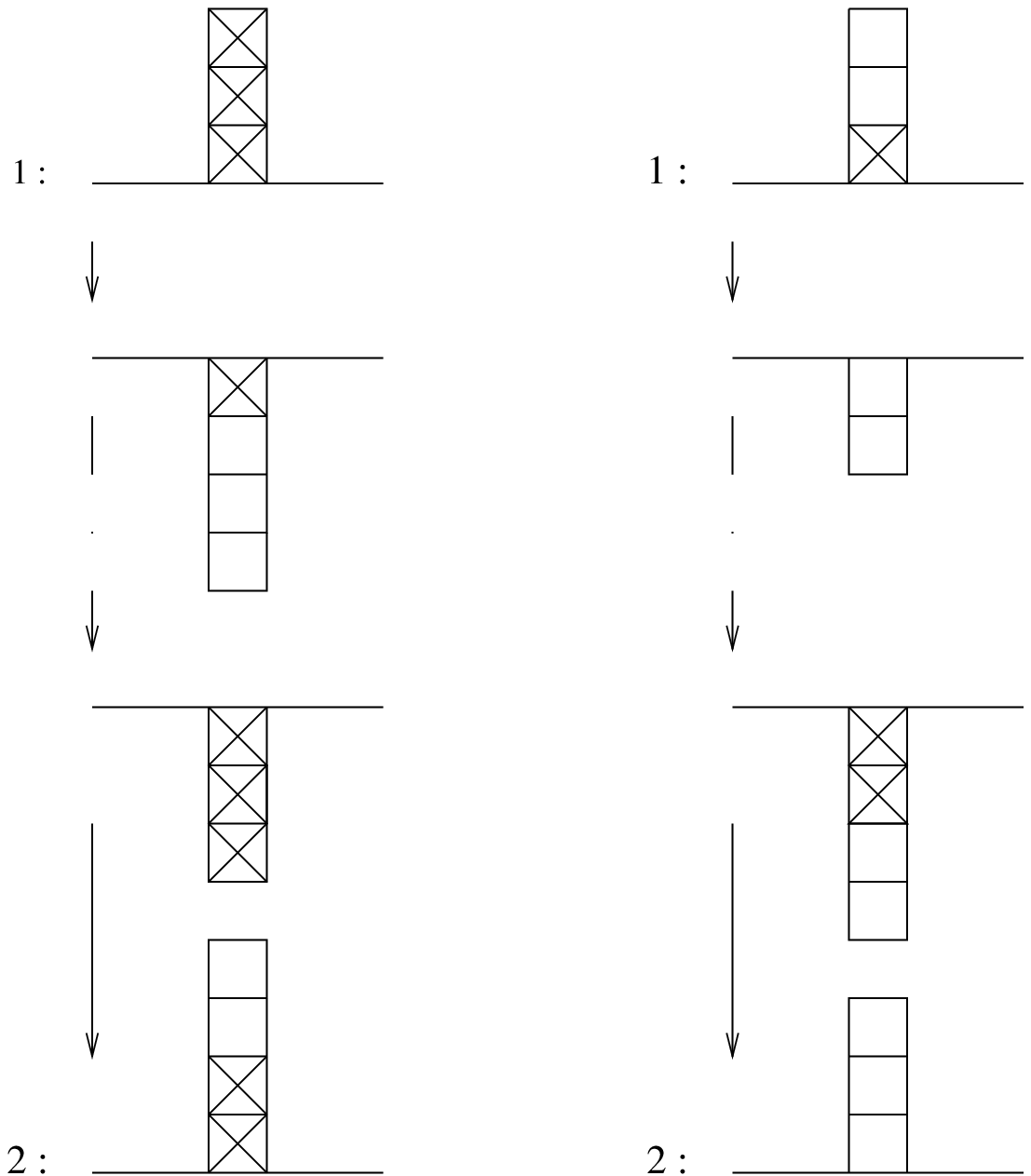}}
\vskip 0.5 cm
\noindent (simple arrows mean single generation, broken arrows mean possibly several generations, but at fixed shape at place $k$).

In view of the characterization of completeness, we observe that we have in both cases: $b_2=b_1+d,\ c_2=c_1+d'$, as if we had not changed the shape (it is easily seen by looking at the $d$ on the left and at the $d'$ on the right).

By Lemma 3, we are now able to remove the condition that the shape does not change at the broken arrows. Indeed, we begin by reasoning about chains as above, then we can ignore the change of shape. By this method we obtain the general result (analogy with a Dick path for which we repeat the removing of sequences $\vee$ and $\wedge$).    
\QED

\vskip 0.2 cm 

{\it Proof of Lemma 3.}
This will be done by induction on $k$.  
\begin{itemize}
\item If $k=1$, the result is obvious. 
\item To prove the result for $k$, we have to show that along the chain between $D$ and $D'$, the shape of the $k$-th column has changed as many times by appearance of a white $x$-column as by appearance of a white $y$-column (i.e., sum of $d$ equal to zero) and as many times
by disappearance of a crossed $x$-column as by disappearance of a crossed $y$-column (i.e., sum of $d'$ equal to zero).

We suppose that our column (assume it is a $y$-column in $D$ and $D'$) changes more times by appearance of a white $x$-column than by appearance of a white $y$-column. Let us observe the subchain on the figure below:

\vskip 0.3 cm
\centerline{
\epsffile{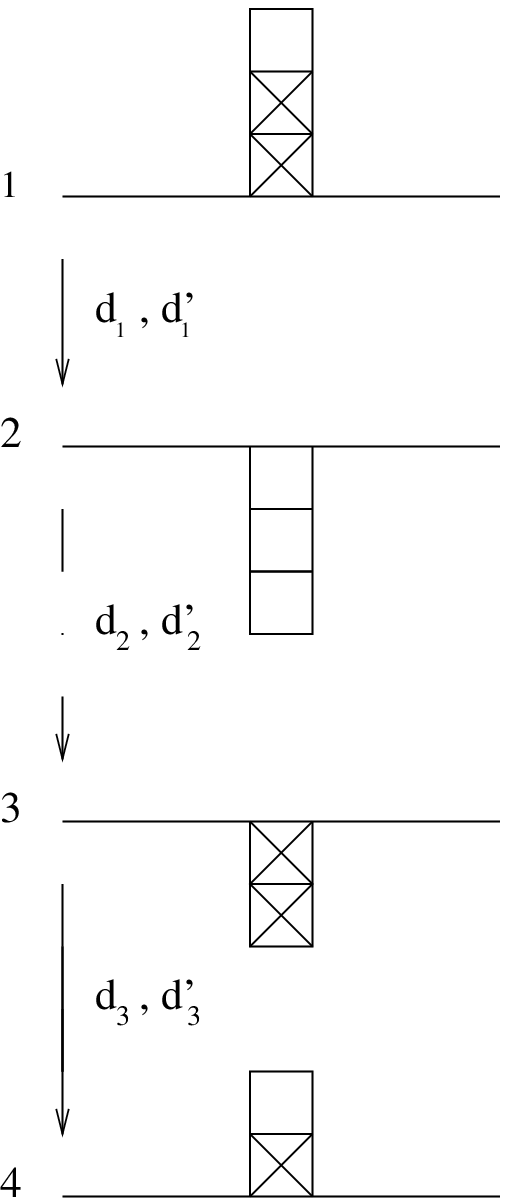}}

Let $h_1$ denote the height of the $y$-column of drawing 1 and $h'_1$ the depth of the first $x$-column on its right. We observe that $b_4=d_3\ge0$ (Case 2 of the Characterization) and that $d_2=h'_1-d_1-d'_1$ since $b_3=0=b_2-d_2$ (Case 1). Thus: $d_1+d_2=h'_1-d'_1$.

We now visualize the changes of shape at place $k$ between $D$ and $D'$ on the following representation.

\vskip 0.3 cm
\centerline{
\epsffile{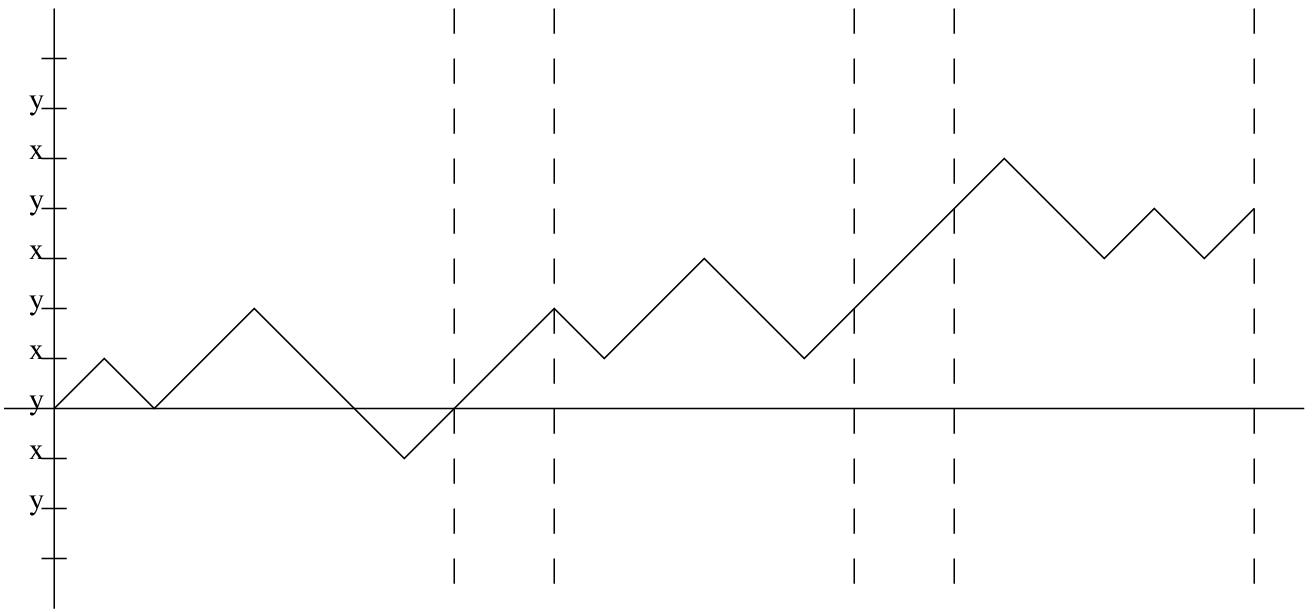}}

The even coordinates correspond to a $y$-column at place $k$, the odd ones to an $x$-column. A north-east line is either the appearance of a white $x$-column or disappearance of a crossed $x$-column (according to odd or even coordinate) and a south-east line is either the appearance of a white $y$-column or disappearance of a crossed $y$-column. The vertical dotted lines are defined as follows. The first is placed at the last point for which the coordinate is equal to zero. Then we have clearly two north-east lines and we put another dotted line. Then we restart with taking coordinate 2 as a new zero for the coordinates. 

Let us suppose that between $D$ and $D'$ there is a single ascent (i.e., a subchain like 1-2-3-4). If we verify that $d_0+d_1+d_2>0$, where $d_0$ is the sum of the $d$ before the ascent, then since $\sum d=0$ between $D$ and $D'$, we have necesseraly some $d<0$ after this sequence, which is impossible without a disappearance of the crossed $y$-column. That is what we wanted to show.

Let us prove that $d_0+d_1+d_2>0$. 

Let $b$ denote the number of white cells at place $k$ of $D$ then $b_1=b+d_0$. Hence: 
$$d_0+d_1+d_2=d_0+h'_1-d'_1=b_1-b+h'_1-d'_1=h_1+h'_1-b.$$
It is easy to check that $h_1+h'_1-b>0$.

It remains to observe that when there are several ascents, the previous reasoning is still true, by looking at the last one. Indeed, it suffices to replace the equality  $b_1=b+d_0$ by $b_1\le b+d_0$ (thanks to what precedes), which keeps the result unchanged.

\end{itemize}

\vskip 0.3 cm

The proof of Lemma 3 is almost complete. It remains to observe that the symmetries between $x$ and $y$ and between crossed and white cells allow us to deal with the other cases. 
\QED

\vskip 0.3 cm
{\bf Lemma 5 :} 

{\it If there is no total completeness along the chain between $D$ and $D'$, then $D\neq D'$ which implies Theorem 3.}

\proof
This is an easy consequence of Lemmas 2 and 4 and Remark 2. It suffices to look at the leftmost place for which the completeness fails: $D'$ has more white cells (and less crosses) than $D$ at this place.
\QED

\subsection{End of the proof}

It is now sufficient to show that there is at least one generation between $D$ and $D'$ that is not complete. We shall in fact show that each generation is not complete.

Let again $D_1=(S_1,T_1)$ and $D_2=(S_2,T_2)$ denote two different drawings, father and son.

If $D_1$ and $D_2$ have the same shape, the result is obvious.

It then remains to study the case where $D_1$ and $D_2$ have different shape. We suppose that completeness holds and reduce it to the absurd.

\vskip 0.2 cm

By looking at the place at most on the left where the shape changes, we can consider only the case where the shape changes at place 1.
The only changes for which the non-completeness is not obvious are the following (remark that here $d=d'=0$):    
\vskip 0.5 cm

\centerline{
\epsffile{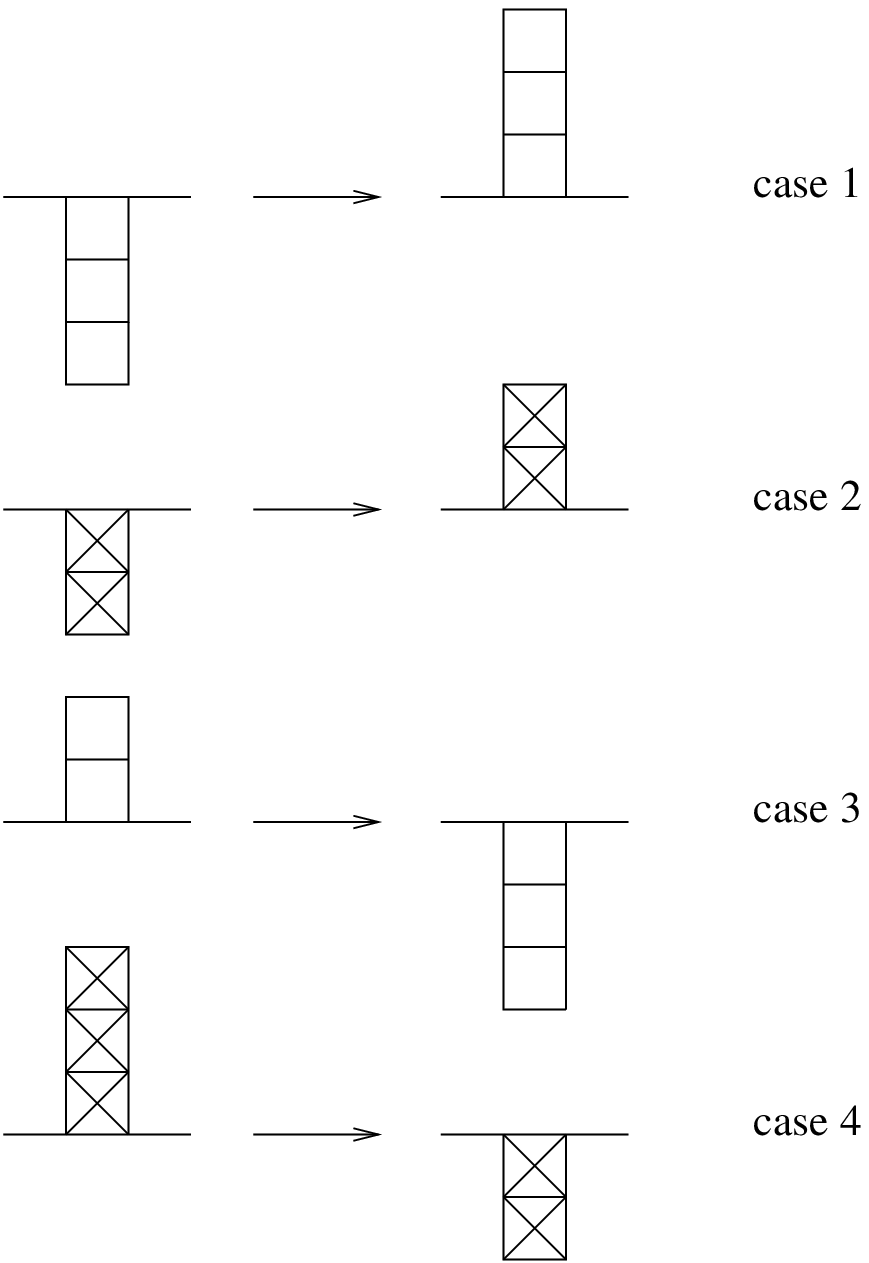}}

\vskip 0.3 cm 
 
The following remark allows us to divide by two the number of cases: 

\vskip 0.2 cm
{\bf Remark 4:}
$D_2$ is a son of $D_1$ if and only if flip($D_1$) is a son of flip($D_2$). This allows us to only consider cases 2 and 4. 

\newpage

\begin{enumerate}
\item Case 2 : 

\centerline{
\epsffile{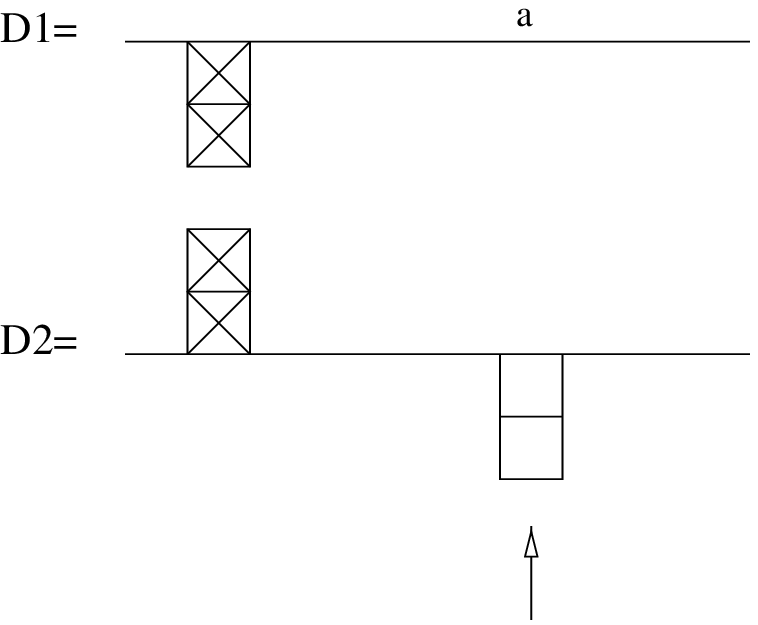}}
If at place ``a''(corresponding to the first plain crossed $x$-column in $D_2$), there is
\begin{itemize}
\item an $x$-column: we first verify that at each place on the left of ``a'' we have $d=0$; then we show that the $x$-column in $D_1$ is smaller than the one in $D_2$, which contradicts $b_2=b_1-d=b_1$.
\item a $y$-column: we first show that in each $x$-column of $T_1+S_2$ there is at least one cell coming from $D_1$ and one coming from $D_2$. This is absurd since there are not enough $x$-columns.
\end{itemize}

\newpage

\item Case 4:

\centerline{
\epsffile{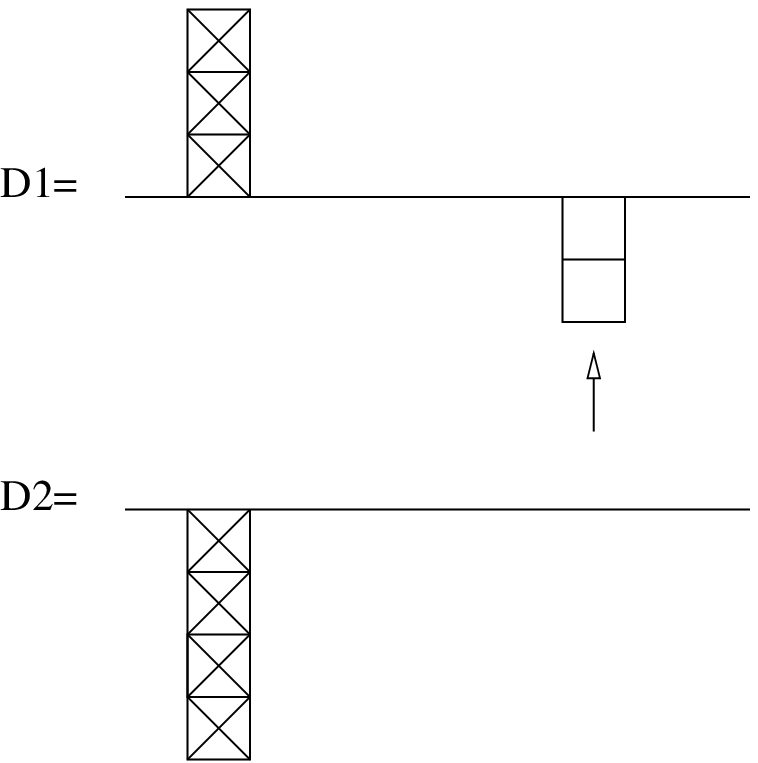}}

In this case, if the first plain $x$-column of $D_1$ is the $l$-th $x$-column of $D_1$, we begin by observing that the $x$-column on its left have at least one white cell, hence have a contribution to $T_1+S_2$. Thus on the left of this place there is already an $x$-column of depth $L-l+1$ (there are at least $l$ $x$-columns in $T_1+S_2$ on the left). This is absurd. 
\end{enumerate}

\section{Elements of $0$ $x$-degree}

\subsection{Description}

Let $\mu=(\mu_1\ge\mu_2\ge\dots\ge\mu_k>0)$ be any partition of $n$.
The goal of this section is to give an explicit basis for $M_{\mu}^0$, which denotes the homogeneous subspace of $M_{\mu}$ of elements of $0$ $x$-degree. We construct this basis with the same visual objects as in the case of hooks. We also obtain a basis for the subspace of $n(\mu)$ $x$-degree which we shall denote by $M_{\mu}^{n(\mu)}$. 

The space $M_{\mu}^0$ has already been studied in [2] and [6]. In particular it is proved that its dimension is $n!/\mu'!$, where $\mu!=\mu_1!\ldots\mu_k!$ and $\mu'$ is the conjugate  of $\mu$. In fact our basis is related to a family introduced in [2]. But we obtain here a direct (and not recursive) method of construction. Moreover we apply the monomial derivatives to $\Delta_{\mu}$ itself and therefore obtain a simple and explicit basis for $M_{\mu}^0$.

We use again the drawings introduced for hook-shaped partitions, here in the case of any partition of $n$. A shape is then made of $n-1$ bars. Each of these bars has $n_x$ $x$-cells and $n_y$ $y$-cells. The set of pairs $(n_x,n_y)$ is the set of biexponents of the partition (the biexponent is omitted). We again put crosses in the shapes and the set of rules for these drawings is the following:
\begin{enumerate}
\item the bars with the same number of $x$-cells are arranged in decreasing height;
\item there are crosses in every $x$-cell;
\item if a bar $B$ is on the left of a bar with more $x$-cells than $B$ and $q$ $y$-cells, then the bar $B$ must have at least $q+1$ $y$-white cells.
\end{enumerate}  

\vskip 0.2 cm
{\bf Remark 5:} By applying flip we obtain a family of drawings with no $x$-crosses.

We now give an example of a drawing:
\vskip 0.5 cm

\centerline{
\epsffile{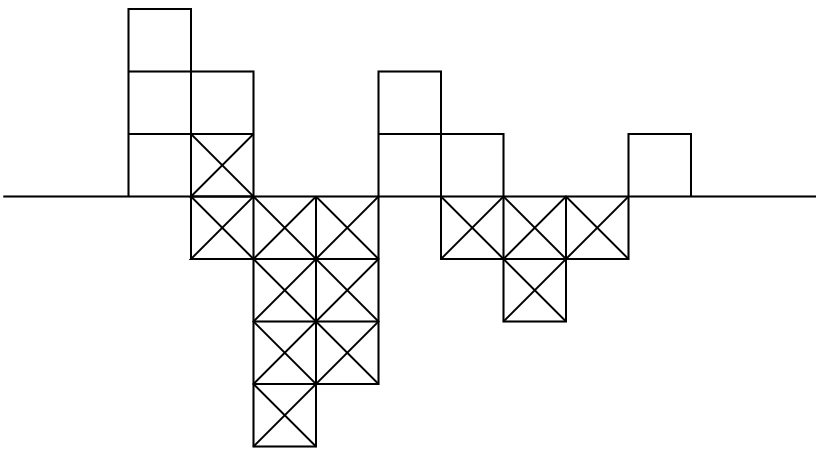}}
\vskip 0.3 cm

\noindent
associated to the partition:
\vskip 0.5 cm

\centerline{
\epsffile{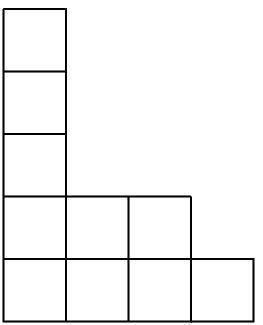}}

\vskip 0.3 cm

\subsection{Enumeration}

We verify that the number of drawings introduced in subsection 1 is $n!/\mu!$. We consider the drawing from the left to the right. The bar we are looking at corresponds to a corner of the Ferrers diagram of the partition from which we have removed the cells corresponding to the bars on the left. 

Number the cells of the partition $\mu$ by writing $i$ in the cell associated to the bar at place $n-i+1$ in the drawing. By the preceding paragraph, this gives a standard tableau.

We now look at the following figures:

\vskip 0.5 cm

\centerline{
\epsffile{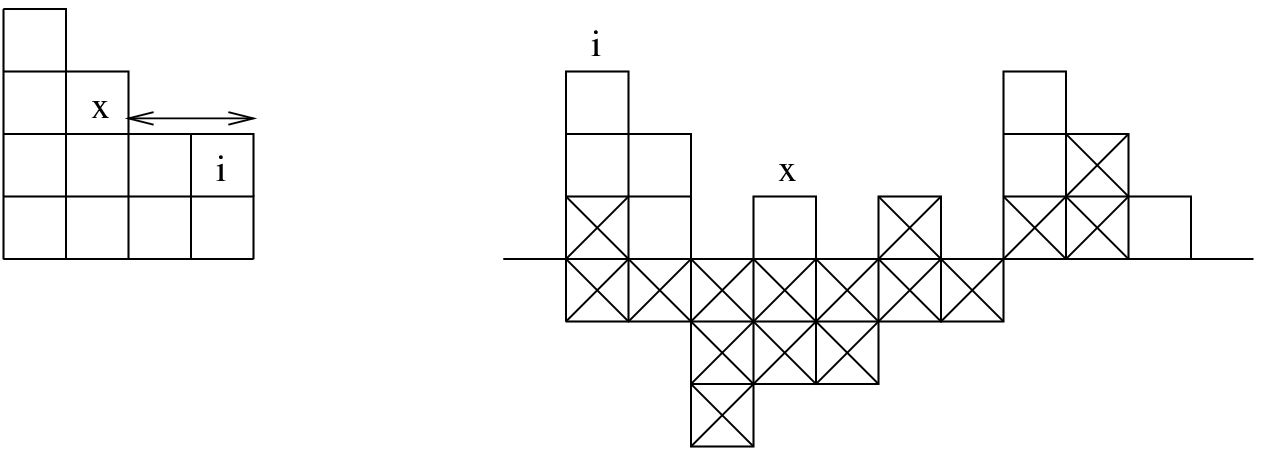}}

We observe that the number of choices for cell $i$ is the length of the arrow that we denote by coarm$^*_i(T^{i+1})$, where $T^{i+1}$ is the (standard) tableau $T$ from which we have removed the cells numbered from $i+1$ to $n$.

We thus obtain that the cardinality is:
$$\sum_{T {\rm standard}} \prod_{i=n}^1 {\rm coarm}^*_i(T^{i+1}). \leqno(1)$$

We show that this number equals $n!/\mu'!$ by induction on $n$. The result is obvious when $n=1$. We write $\mu'=(c_1^{\alpha_1},\ldots,c_h^{\alpha_h})$, where the $c_j$'s are  the height of columns of $\mu$ and $\alpha_j$ their multiplicities. In particular, $\mu$ has $h$ corners, $\mu'!=\prod_{i=j}^{h}(c_j!)^{\alpha_j}$, $n=\sum_{j=1}^h \alpha_j c_j$ and $\alpha_j$ is the contribution of corner $j$ in the product of (1). We then rewrite this formula as:
$$\sum_{j=1}^h \alpha_j.\sum_{T'} \prod_{i=n-1}^1 {\rm coarm}^*_i(T'^{i+1})$$ 
where $T'$ varies amongst every standard tableaux of the Ferrers diagram from which we have removed its $j$-th corner (let $\mu^j$ denote the corresponding partition). We are now able to conclude, since $\mu'^j!=\mu'!/c_j$:
$$\sum_{j=1}^h \alpha_j {\frac {(n-1)!} {\mu'^j!}}={\frac {(n-1)!} {\mu'!}}\sum_{j=1}^h \alpha_j.c_j={\frac {n!} {\mu'!}}.$$

\subsection{Independence and conclusion}

As in the case of the hook-shaped partitions, we denote by $S$ (respectively $T$) the diagram consisting only of the crosses (respectively of the white cells) of a given drawing. For example in the case of the drawing of subsection 5.1, we have:

\vskip 0.5 cm

\centerline{
\epsffile{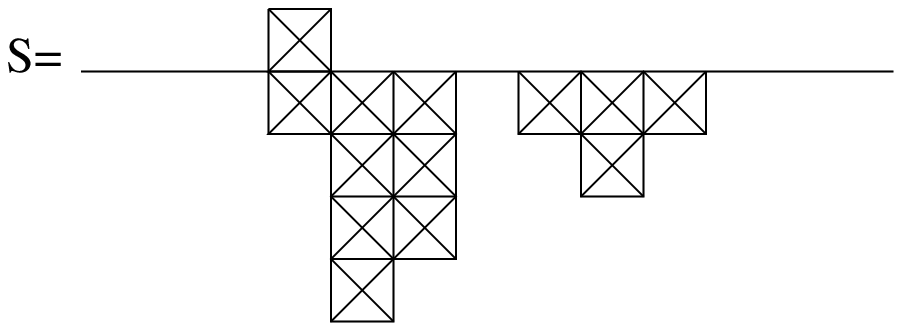}}

\vskip 0.5 cm

\centerline{
\epsffile{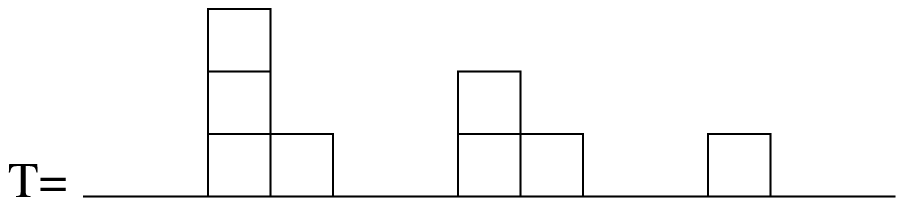}}

\vskip 0.3 cm

Let now $\partial_S$ and $\partial_T$ denote the derivative operators associated to $S$ and $T$. Let also ${\cal S}$ and ${\cal T}$ denote the set of all $S$'s and $T$'s constructed by this way. We also associate to $S$ and $T$ a monomial by the same way. For example, in the case of the previous figures, we have: $M_S=x_2y_2x_3^4x_4^3x_6x_7^2x_8$ and $M_T=y_1^3y_2y_5^2y_6y_9$.

\vskip 0.3 cm
{\bf Theorem 4:}

{\it The set $\{\partial_S\Delta_{\mu}\}_{S\in{\cal S}}$ is linearly independent and hence is a basis of $M_{\mu}^0$.

The set $\{\partial_T\Delta_{\mu}\}_{T\in{\cal T}}$ is linearly independent and hence is a basis of $M_{\mu}^{n(\mu)}$.}

\vskip 0.3 cm
The theorem is a consequence of the two following lemmas.

\vskip 0.3 cm
{\bf Lemma 6:}

{\it We can reconstruct the drawing from $S$ or $T$.}

\newpage

{\bf Lemma 7:}

{\it For the lexicographic order ($x_1<x_2<\dots<x_n<y_1<\dots<y_n$) $M_T$ is the minimal monomial for $\partial_S\Delta_{\mu}$ and $M_S$ for $\partial_T\Delta_{\mu}$.}

\vskip 0.3 cm

\proof
The proof of Lemma 6 is easy: we reconstruct the drawing from the left to the right, as in the case of hooks, thanks to the rules.

The proof of Lemma 7 requires attention only in the case of $T$, so we develop this point. Once the crossed cells have been fixed, we have to show that the white cells are at most on the left. It suffices in fact to show that the $x$-white cells can not be moved to the left. We show it by looking at the drawing from the left to the right. Let $k$ and $l$ denote the number of $x$-cells and of $y$-crossed cells at place $p$. We have to prove that a bar with $l'>l$ $x$-cells and $k'\ge k$ $y$-crossed cells is forbidden at place $p$. If the couple $(k',l')$ is not a biexponent of the partition or if it is present on the left, we are done. To conclude we observe that this couple can not be a biexponent of the partition appearing on the right of the initial drawing. Indeed, because of the rules we should have: $k>k'$. 
\QED

\vskip 0.2 cm
{\bf Remark 6:} It is possible to show that our family of monomials $\{M_S\}_{S\in {\cal S}}$ is equal to the family $B_{\mu}$ of [2], section 4. But whereas $B_{\mu}$ was constructed recursively, our construction is direct. Moreover we apply it directly to $\Delta_{\mu}$ and obtain simple and explicit bases for $M_{\mu}^0$ and $M_{\mu}^{n(\mu)}$, whereas N. Bergeron and A. Garsia were dealing in [2] with linear translates of Garnir polynomials.

\section{References}
\begin{enumerate}

\item E. Allen, {\it The decomposition of a bigraded left regular representation of the diagonal action of $S_n$}, J. Comb. Theory A, {\bf 71} (1995), 97-111.

\item N. Bergeron and A. M. Garsia, {\it On certain spaces of harmonic polynomials}, Contemporary Mathematics, {\bf 138} (1992), 51-86.

\item Louis Comtet, {\it Analyse Combinatoire}, Presses Universitaires de France, Paris, 1970.

\item A. M. Garsia and M. Haiman, {\it Orbit harmonics and graded representation}, in ``Laboratoire de combinatoire et d'informatique math\'ematique, UQAM collection'' (S. Brlek, Ed), to appear.

\item A. M. Garsia and M. Haiman, {\it A graded representation model for Macdonald's polynomials}, Proc. Natl. Acad. Sci., {\bf 90} (1993), 3607-3610.

\item A. M. Garsia and M. Haiman, {\it Some natural bigraded $S_n$-modules and $q,t$-Kostka coefficients}, Elec. J. of Comb. 3 (no. 2) (1996), R24. 

\item A. M. Garsia and J. Remmel, {\it Plethystic formulas and positivity for $q,t$-Kostka polynomials}, In Mathematical Essays in Honor of Gian-Carlo Rota (Cambridge, MA, 1996), Birkh\"auser Boston, Boston, MA (1998), 245-262.

\item A. M. Garsia and G. Tesler, {\it Plethystic formulas for Macdonald $q,t$-Kostka coefficients}, Advances in Math., {\bf 123} (1996), 144-222.

\item M. Haiman, {\it Macdonald polynomials and geometry}, preprint.

\item A. N. Kirillov and M. Noumi, {\it Affine Hecke algebras and raising operators for Macdonald polynomials}, Duke Math. J., {\bf 93} (1998), 1-39.

\item F. Knop, {\it Integrality of two variable Kostka functions}, J. Reine Angew. Math., {\bf 482} (1997), 177-189.

\item I. G. Macdonald, {\it A new class of symmetric functions}, Actes du $20^e$ S\'eminaire Lotharingien, Publ. I.R.M.A. Strasbourg (1988), 131-171.

\item E. Reiner, {\it A Proof of the $n!$ Conjecture for Generalized Hooks}, J. Comb. Theory A, {\bf 75} (1996), 1-22.

\item S. Sahi, {\it Interpolation, integrality, and a generalization of Macdonald's polynomials}, Internat. Math. Res Notices, {\bf 10} (1996), 457-471.

\end{enumerate}

\noindent
{\large{\bf  Acknowledgements}}

\medskip

The author would like to express all his gratitude to the referees who have made significant efforts to improve this paper by their valuable advices and suggestions.

\end{document}